\documentclass{notices}
\usepackage{amsmath,amsfonts,amsthm,amscd,amssymb,graphicx,color}

\newcommand{\rmi}{\mathrm{i}}

\title{Spectral stability and spatial dynamics in partial differential equations
}

\author{
Margaret Beck\footnote{
    Margaret Beck is a professor of mathematics at Boston University. Her email address is mabeck@bu.edu.
    }
}

\begin{document}

\maketitle

\section*{} % Effectively an introduction, but not labeled as such.

This article is focused on two related topics within the study of partial differential equations (PDEs) that illustrate a beautiful connection between dynamics, topology, and analysis: \emph{stability} and \emph{spatial dynamics}\footnote{{\color{black}The connection between these two concepts was also described in the talk entitled ``Stability for PDEs, the Maslov Index, and Spatial Dynamics," which the author gave at MSRI in 2018. That talk can be accessed via \url{https://www.msri.org/workshops/871/schedules/24652}}}. The first is a property of solutions that describes the extent to which they can be expected to persist, and hence be observed, over long time scales. The second is a perspective that has been used to study various properties, such as stability, of nonlinear waves and coherent structures, the term often used to describe the solutions of interest in the class of PDEs that will be considered here. 

To fix ideas, let's focus on systems of reaction-diffusion equations, 
\begin{equation}\label{E:main-rd}
u_t = \Delta u + f(u),
\end{equation}
where $u: \Omega \times [0, \infty) \to \mathbb{R}^n$, $\Omega \subset \mathbb{R}^d$, $f: \mathbb{R}^n \to \mathbb{R}^n$, $\Delta = \nabla \cdot \nabla = \partial_{x_1}^2 + \dots + \partial_{x_d}^2$, and there are accompanying initial conditions and possibly also boundary conditions for $\partial \Omega$, which for the moment I will leave unspecified. I will assume $f$ and $\partial \Omega$ are smooth.  

Reaction-diffusion equations are a class of parabolic PDEs for which it is interesting to study the dynamics specifically because well-posedness is known: under reasonably mild assumptions, unique solutions exist and depend smoothly on the initial data and the function $f$. This means that one can focus on the resulting behavior of solutions as time evolves, and in many cases obtain quite detailed information. They are also {\color{black} relevant} because they appear in a wide variety of applications, for example in chemistry, biology, and ecology, which means that not only are there specific models in which to test the theory, but there are also important open questions originating in other sciences that can point to interesting new mathematical directions. 

It is worth noting that many of these properties that have just been described are also present in other types of PDEs, such as the nonlinear Schr\"odinger equation and the Korteweg-de Vries equation, both dispersive evolution equations, and so much of what will be discussed below can be applied not just to reaction-diffusion equations but also more broadly. See \cite{ChardardDiasBridges11} for a variety of examples related to the context of this article.

\section*{Stability}

In order to describe the dynamics of the PDE one often begins by identifying specific solutions, such as stationary or time-periodic patterns, and then seeking to understand the extent to which such solutions will be observed in the long-time dynamics. Within this context, one might ask about two types of stability. The first is related to robustness of the solution to perturbation in the system parameters, or in other words to perturbations within the PDE itself. This type of stability is referred to as structural stability, and it typically falls within the sub-field of bifurcation theory. The second type of stability, and the one that is a focus of this article, is stability in time, or dynamic stability: can one expect to observe this solution in the dynamics of a fixed PDE as time evolves? This has to do with robustness of the solution to perturbations in the initial condition, or to perturbations in the current state of the system. In this sense, stable solutions attract (or at least do not repel) nearby data. Unstable states repel (at least some) nearby data, which will be driven away to {\color{black}some structure} that is dynamically attracting. Structural and dynamic stability are of course connected; one could for example ask how dynamic stability is affected by changes in system parameters. But for the remainder of this article, stability will always refer to stability in time. 

Let's suppose that we are given a stationary solution of \eqref{E:main-rd}, $\varphi(x)$, so that
\begin{equation}\label{E:stationary}
0 = \Delta \varphi + f(\varphi),
\end{equation}
and we want to investigate its stability. We can write the solution to \eqref{E:main-rd} as $u(x,t) = \varphi(x) + v(x,t)$ and derive an evolution equation for the perturbation $v$:
\[
v_t = \underbrace{\Delta v + df(\varphi) v}_{=: \mathcal{L} v} + \underbrace{[f(\varphi + v) - f(\varphi) -  df(\varphi)v]}_{=: \mathcal{N}(v)}.
\]
If $v(x,0)$ is small in some appropriate sense (so we are focusing on local, rather than global, stability), will the perturbation decay to zero, or at least remain small, for all $t \geq 0$? 

Because the perturbation $v$ is small, at least initially, one could expect the linear term $\mathcal{L} v$ to dominate the nonlinear one $\mathcal{N}(v)$ in determining the dynamics, simply because $|v|^p < |v|$ if $p > 1$ and $|v| <1$. Thus, one could focus initially on the linear dynamics, in which case the spectrum of $\mathcal{L}$ plays a key role. This relies on the fact that the linear operator is nice: it generates an analytic semigroup, and so there is a clear connection between spectrum and dynamics. Unstable (positive real part) spectrum leads to exponential growth, stable (negative real part) spectrum leads to exponential decay, and if there is spectrum on the imaginary axis then one must take the nonlinearity into account. 

Here the focus will be on detecting spectral instabilities. The spectrum of $\mathcal{L}$ can be divided into two parts: the essential spectrum and the point spectrum, or eigenvalues. At the moment the details of this decomposition are not so important; what is important is the fact that the essential spectrum is relatively easy to compute, whereas the point spectrum is typically difficult to compute. Thus, if one calculates the essential spectrum and it lies in the right half plane, then an instability has been detected. The more interesting {\color{black} case is therefore} when the essential spectrum is stable, and one needs to understand the point spectrum. Thus, the question of detecting an instability is reduced to determining whether or not there are any eigenvalues of the linearized operator that have positive real part.  

The simplest case is a scalar equation in one space dimension: $n = d = 1$. If $\Omega = (a,b)$ and we consider zero Dirichlet boundary conditions, then we are in the classical setting of a Sturm-Liouville eigenvalue problem:
\begin{gather*}
\lambda v = v_{xx} + df(\varphi(x))v, \qquad x \in (a, b) \\
v(a) = v(b) = 0. 
\end{gather*} 
Note the linear operator is self-adjoint, so the spectrum is real\footnote{{\color{black} On the bounded domain considered here, one could attribute the realness of the spectrum to the fact that the operator is second-order and scalar, since any second-order scalar operator can be put into self-adjoint form by means of an appropriate integrating factor. Later, however, we will consider operators on the entire real line that act on vector-valued functions, in which case the realness of the spectrum will result from the self-adjointness of the operator.}}. {\color{black} Consider the Pr\"ufer coordinates
\[
v = r \sin \theta, \qquad v_x = r \cos \theta,
\]
which {\color{black} in this setting are essentially just polar coordinates in the phase plane}. By differentiating the relations $r^2 = v^2 + v_x^2$ and $\tan \theta = v/v_x$ and solving for $r_x$ and $\theta_x$}, we find the dynamics of $r$ and $\theta$ to be governed by
\begin{eqnarray*}
r_x &=& r(1 + \lambda - df(\varphi(x))) \cos \theta \sin \theta, \\ 
\theta_x &=& \cos^2 \theta + (df(\varphi(x)) - \lambda) \sin^2 \theta.
\end{eqnarray*}
One can now make three key observations: the dynamics for $\theta$ have decoupled from those for $r$; the set $\{ r = 0\}$ is invariant; and therefore a solution that is not identically zero can satisfy the boundary condition only if $\theta(a; \lambda), \theta(b; \lambda) \in \{ j \pi\}_{j \in \mathbb{Z}}$. Thus, the second order eigenvalue problem has been reduced to the study of the first order equation for $\theta$: if for a given $\lambda$ there exists a solution $\theta$ satisfying the boundary condition, then $\lambda$ is an eigenvalue of $\mathcal{L}$.  

Let's shift our perspective slightly and, rather than thinking of $x$ as a spatial variable, let us view it as a time-like variable. (This is an example of spatial dynamics.) If $ \theta(a; \lambda) \notin \{ j \pi\}_{j \in \mathbb{Z}}$, then $\theta$ cannot be an eigenfunction; therefore to determine if $\lambda$ is an eigenvalue, by periodicity we can assume $\theta(a; \lambda) = 0$. Because of the structure of the equation, for $\lambda$ large and negative we expect $\theta$ to oscillate and to find eigenvalues. Suppose we have found one, and we label it $\lambda_k$ to indicate $\theta(b; \lambda_k) = (k+1)\pi$. If we continuously increase $\lambda$, we continuously decrease $\theta(b; \lambda)$, and the next eigenvalue occurs when we reach the point where $\theta(b; \lambda_{k-1}) = k \pi$. Expanding on this argument, one can prove there is a sequence of simple eigenvalues $\lambda_0 > \lambda_1 > \dots$ and corresponding sequence of solutions $\theta$ such that $\theta(b; \lambda_k) = (k+1)\pi$. This in turn implies that the corresponding eigenfunction $v(x; \lambda_k)$ has exactly $k$ simple zeros in the interval $(a,b)$.

\begin{figure}
\centering
\includegraphics[width=0.45\textwidth]{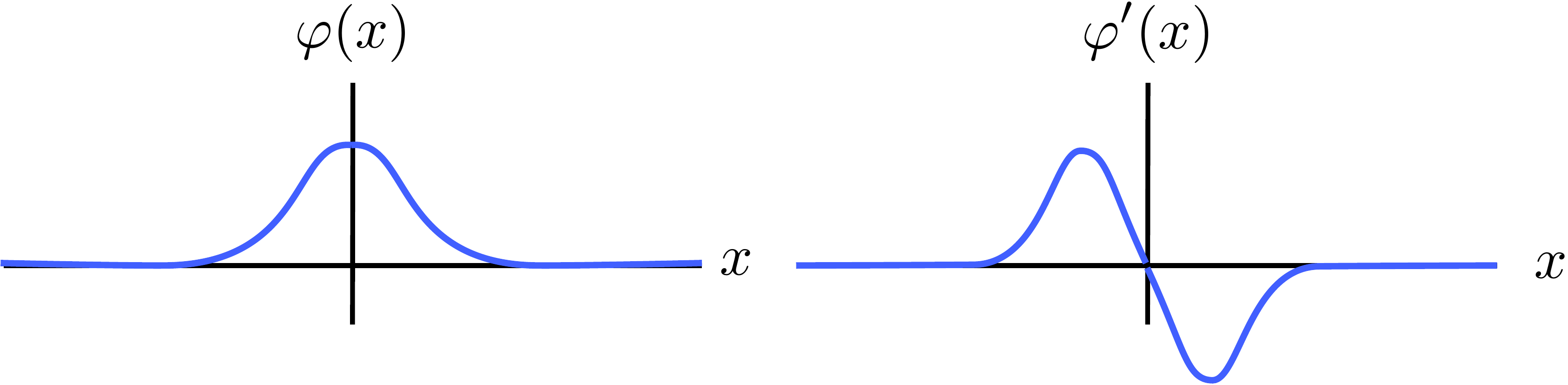}
\caption{A pulse and its derivative.}
\label{F:pulse}
\end{figure}

From the perspective of stability, this is an extremely powerful result. This is classically illustrated by considering a scalar reaction-diffusion equation on the entire real line that has a pulse as a stationary solution; see Figure \ref{F:pulse}. {\color{black} This is a natural example to consider for at least two reasons. First, in the context of applications reaction-diffusion equations are often posed on the entire real line so as to avoid any potential complications arising from the boundary while still capturing the experimentally observed behavior. Second, pulses are among the simplest and most common type of coherent structures found in such models. The relevant elements of the above theory remain when we replace the interval $(a,b)$ with the real line $\mathbb{R}$, as long as we work in an appropriate function space, such as $L^2(\mathbb{R})$.} Because $\varphi$ satisfies \eqref{E:stationary}, if we take an $x$-derivative of this equation we find that $0 = \mathcal{L} \varphi_x$, and so $\varphi_x$ is an eigenfunction of $\mathcal{L}$ with eigenvalue zero. As illustrated in Figure \ref{F:pulse}, $\varphi_x$ has exactly one zero. This implies that $0 = \lambda_1$, and so there must be a positive eigenvalue, $\lambda_0 > 0$. As a result, any stationary pulse solution of a scalar reaction-diffusion equation on the real line must be unstable. The details of the function $f$ are not relevant, other than that the resulting equation has a pulse solution, nor are the details of $\varphi$, other than that it is a pulse (or more generally has at least one local extrema). A complementary result holds if $\varphi$ is a monotonic front, in which case $\varphi_x$ has no zeros, and so the largest eigenvalue is  zero: $\lambda_0 = 0$.

In this example, the zeros of the eigenfunction are being used as a proxy for the eigenvalues. This suggests the alternative perspective of conjugate points, which can be described as follows. Above, the domain $(a,b)$ was kept fixed, $\lambda$ was allowed to vary, and the values of $\lambda$ where the solution satisfied the boundary condition were recorded. Instead, let's fix $\lambda$ and allow the domain to vary: $x \in (a, s)$ with $s \in [a, b]$. {\color{black} The number $s$ is defined to be a conjugate point for $\lambda$ if $\lambda$ is an eigenvalue of the Dirichlet problem posed on the domain $[a, s]$.} We can play a similar game if we fix $\lambda = \lambda_k$. We therefore know that if $s = b$, then $\theta(b; \lambda_k) = (k+1)\pi$. We can now continuously decrease $s$ from $b$, so that $\theta$ has less time to oscillate (that's the spatial dynamics perspective again), and record the values $s_j$ where $\theta(s_j; \lambda_k) = (j+1)\pi$. In this way, we get a sequence of conjugate points $s_k = b > s_{k-1} > s_{k-2} > \dots > s_0 > a$ that are in one-to-one correspondence with the eigenvalues that are strictly bigger than $\lambda_k$. 

\begin{figure}
\centering
\includegraphics[width=0.45\textwidth]{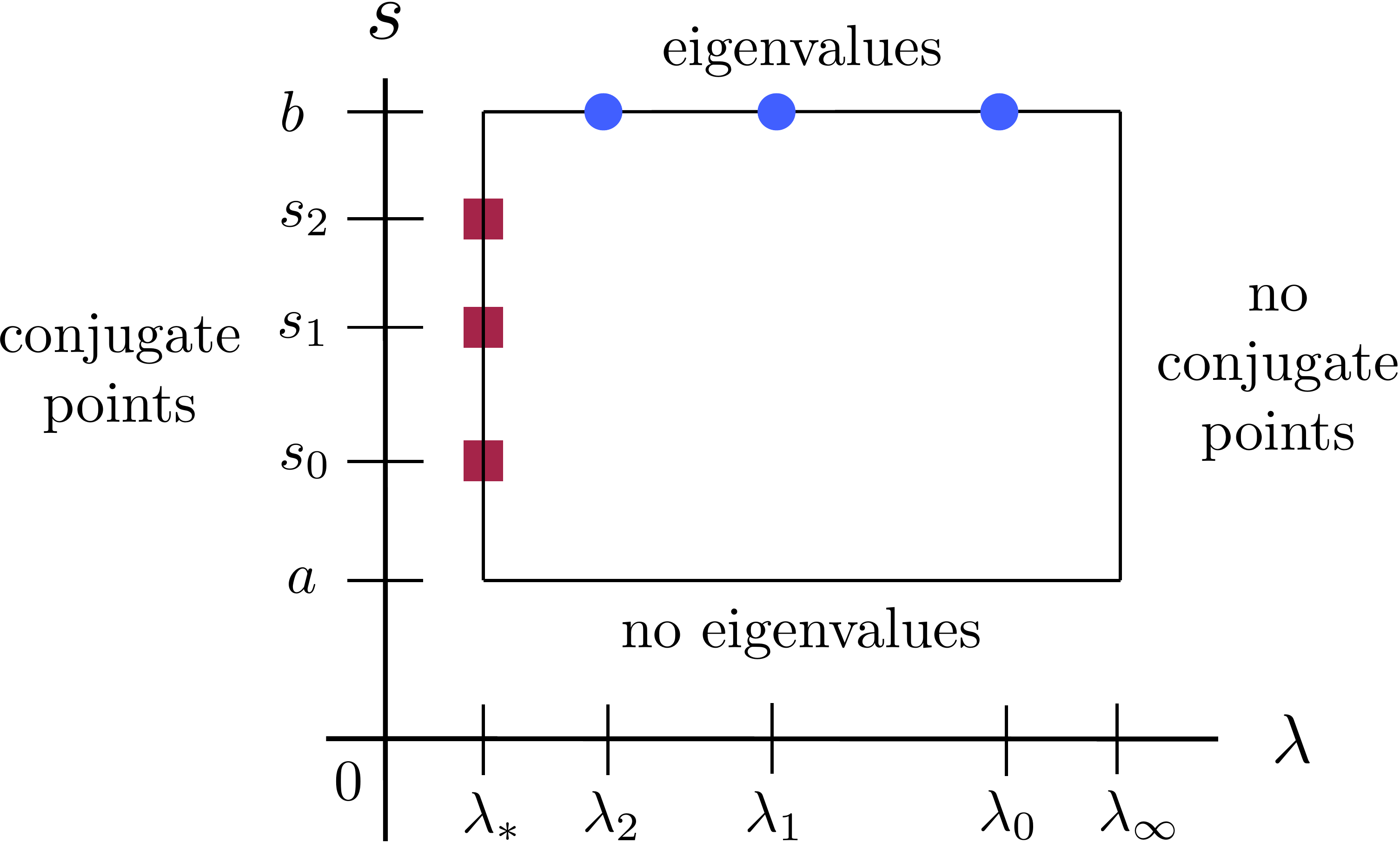}
\caption{The square illustrating that the number of conjugate points for $\lambda = \lambda_*$ is equal to the number of eigenvalues $\lambda > \lambda_*$.}
\label{F:square}
\end{figure}

This result is illustrated using the ``square" depicted in Figure \ref{F:square}. To complete the picture, one needs to show that for $\lambda = \lambda_\infty$ sufficiently large there are no conjugate points, and note that for $s = a$ there are no eigenvalues simply because there are no dynamics. To detect instabilities, one can fix $\lambda_* = 0$, and then the number of conjugate points must be equal to the number of unstable eigenvalues. In the example above regarding pulse instability, by counting zeros of $\varphi_x$ we were effectively counting conjugate points to prove the existence of an unstable eigenvalue. This is a simple case of what's often called the Morse Index Theorem, and it goes back to the work of Morse \cite{Morse96}, Bott \cite{Bott56}, and others. 

The idea of counting unstable eigenvalues by instead counting conjugate points seems nice, but it appears to be restricted to the scalar case, where we can use polar coordinates to define the angle $\theta$. However, Arnol'd \cite{Arnold85,Arnold67} realized that a generalization of this angle to the system case ($n > 1$) was possible using the Maslov Index, and that this enabled the study of the associated oscillations; his ideas were then utilized in \cite{Jones88} to prove instability of a standing wave in a nonlinear Schr\"odinger-type equation. This latter paper was the catalyst for the ideas which we now describe.

Let's return to equation \eqref{E:main-rd}, but retain the restriction to one space dimension: $d=1$. To most directly utilize the Maslov index, we'll assume the nonlinearity is a gradient, $f = \nabla G$ for some $G: \mathbb{R}^n \to \mathbb{R}$. The eigenvalue problem then becomes
\[
\lambda v = v_{xx} + \nabla^2G(\varphi(x)) v = \mathcal{L} v, \qquad x \in \mathbb{R},
\]
where now $\Omega = \mathbb{R}$ and it is required that $v \in L^2(\mathbb{R}; \mathbb{R}^n)$, in lieu of specifying boundary conditions. Note that the linear operator is again self adjoint, so $\lambda \in \mathbb{R}$. To fix ideas, let's again suppose $\varphi$ is a pulse, meaning that $\lim_{x \to \pm \infty}\varphi(x) = \varphi_\infty$ for some $\varphi_\infty \in \mathbb{R}^n$. As mentioned above, the most interesting case is to assume the essential spectrum of $\mathcal{L}$ is stable, so we can focus on detecting unstable eigenvalues. It turns out this is equivalent to assuming that $\nabla^2G(\varphi_\infty)$ is a negative matrix; this will be utilized below. This second-order eigenvalue problem can again be written as a first order system, now via
\begin{equation}\label{E:eval-sys}
\frac{d}{dx} \begin{pmatrix} v \\ w \end{pmatrix} = \underbrace{\begin{pmatrix} 0 & -1 \\ 1 & 0 \end{pmatrix}}_{=:J} \underbrace{\begin{pmatrix} \lambda - \nabla^2G(\varphi(x)) & 0 \\ 0 & -I \end{pmatrix}}_{=: \mathcal{B}(x; \lambda)} \begin{pmatrix} v \\ w \end{pmatrix}. 
\end{equation}
There's that spatial dynamics perspective again. 

To understand how to associate an angle with this first-order eigenvalue problem, let's step back and discuss the Maslov Index. {\color{black} An accessible explanation of the topics we are about to describe can be found in \cite{HowardLatushkinSukhtayev17}.} To begin, consider the symplectic form $\omega(U, V) := \langle U, JV \rangle_{\mathbb{R}^{2n}}$, where $J$ is defined in \eqref{E:eval-sys} and $\langle \cdot, \cdot \rangle_{\mathbb{R}^{2n}}$ is the usual inner product in $\mathbb{R}^{2n}$. The associated Lagrangian-Grassmanian is the set of all $n$-dimensional subspaces of $\mathbb{R}^{2n}$ on which the symplectic form vanishes:
\[
\Lambda(n) = \{ \ell \subset \mathbb{R}^{2n}: \mathrm{dim}(\ell) = n, \mbox{ } \omega|_{\ell \times \ell} = 0 \}.
\]
Each Lagrangian plane has an associated frame matrix, defined in terms of square matrices $A, B \in \mathbb{R}^n$ such that 
\[
\ell = \left \{ \begin{pmatrix} A \\ B \end{pmatrix} u: u \in \mathbb{R}^n\right\}.
\]
{\color{black} The plane is just the column space of the frame matrix. In fact, the above frame matrix is not unique, and each plane corresponds to an equivalence class of frame matrices.} Suppose we have a path of Lagrangian subspaces, $\ell(t)$ for $t \in (a, b)$, and we are interested in intersections of this path with a fixed reference Lagrangian plane, say the Dirichlet plane: $\mathcal{D} = \{ (0, v) \in \mathbb{R}^{2n}: v \in \mathbb{R}^n\}$. (This is {\color{black} analogous to} looking for conjugate points.) Associate the path $\ell(t)$ with frame matrices $A(t), B(t)$. Arnol'd showed there is a well-defined angle $\theta(t)$ such that
\begin{equation}\label{E:def-theta}
e^{\rmi \theta(t)} = \mathrm{det}[\underbrace{(A(t) - \rmi B(t))(A(t) + \rmi B(t))^{-1}}_{=:W(t)}].
\end{equation}
The reason this works is that the Lagrangian structure of $\ell$ forces $W$ to be unitary, so its spectrum lies on the unit circle. Moreover, it can be shown that 
\[
\mathrm{dim }[\mathrm{ ker}(W(t) + I)] = \mathrm{dim}(\ell(t) \cap \mathcal{D}).
\]
{\color{black} Note that the quantity on the left hand side refers to the complex dimension of the complex vector space $\mathrm{ ker}(W(t) + I) \subset \mathbb{C}^n$, whereas the quantity on the right hand side refers to the real dimension of the real subspace of $\ell(t) \cap \mathcal{D} \subset \mathbb{R}^{2n}$.}
To write down the definition of the Maslov index in full detail would be quite lengthy; here the key fact is that the Maslov index counts, with multiplicity and direction, the number of times an eigenvalue of $W(t)$ crosses through $-1$. Hence, it is also counting intersections of the path $\ell(t)$ with the reference plane $\mathcal{D}$. In this sense, the Maslov index counts conjugate points. 
The Maslov index is related to the fact that the fundamental group of the Lagrangian-Grassmanian is the integers; if $\ell(t)$ is a loop, its Maslov index is its equivalence class in the fundamental group \cite{Arnold67}. %\cite{HowardLatushkinSukhtayev17}. %\cite{Furutani04}

\begin{figure}
\centering
\includegraphics[width=0.5\textwidth]{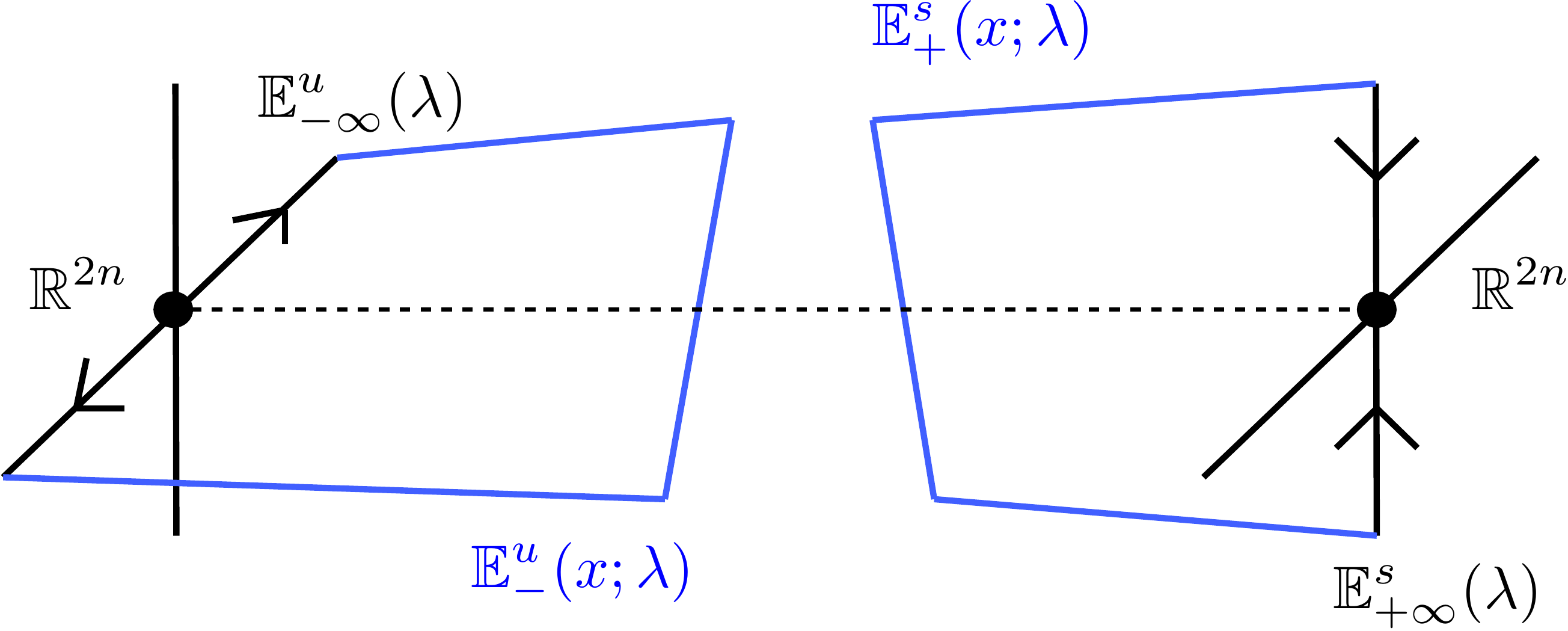}
\caption{The subspaces of decaying solutions.}
\label{F:intersecting-subspaces}
\end{figure}

Let's return now to our eigenvalue problem \eqref{E:eval-sys}. Our assumption that the essential spectrum is stable, $\nabla^2 G(\varphi_\infty) < 0$, implies that the asymptotic matrices $\lim_{x\to\pm\infty}J\mathcal{B}(x; \lambda)$ are both hyperbolic, with stable and unstable subspaces of dimension $n$. If we let $\mathbb{E}^u_-(x; \lambda)$ and $\mathbb{E}^s_+(x; \lambda)$ denote the subspaces of solutions that are asymptotic to the unstable eigenspace at $-\infty$ and the stable subspace at $+\infty$, respectively, then in order to have an eigenfunction $v \in L^2$ we must have $(v, w)(x; \lambda) \in \mathbb{E}^u_-(x; \lambda) \cap \mathbb{E}^s_+(x; \lambda)$; otherwise, the solution would be growing exponentially fast in forward or backward time, thus preventing $v$ from being square integrable. See Figure \ref{F:intersecting-subspaces}.

Studying the intersection of these subspaces leads to the now standard theory behind the Evans function \cite{Sandstede02}. So far we have made no reference to any Lagrangian structure. It turns out that our assumption that $f = \nabla G$ implies that in fact both $\mathbb{E}^u_-(x; \lambda)$ and $\mathbb{E}^s_+(x; \lambda)$ are paths of Lagrangian subspaces. With this additional structure, we can adopt a different perspective and look for conjugate points: given $\ell(x; \lambda) := \mathbb{E}^u_-(x; \lambda) \in \Lambda(n)$, we define a conjugate point to be a value of $x$ such that $\ell(x; \lambda) \cap \mathcal{D} \neq \{0\}$.

Using this framework, in \cite{BeckCoxJones18} it was shown that the square depicted in Figure \ref{F:square}, suitably adapted to reflect the fact that the spatial domain is now all of $\mathbb{R}$, holds for the eigenvalue problem \eqref{E:eval-sys}. {\color{black} This relies on the homotopy invariance of the Maslov index and the fact that the boundary of the square maps to a null-homotopic curve in the Lagrangian-Grassmanian.} Thus, one can count unstable eigenvalues by instead counting conjugate points. Furthermore, this result was used to prove that, in equations of the form \eqref{E:main-rd} with $\Omega = \mathbb{R}$ and $f = \nabla G$, any {\color{black} generic} pulse solution must necessarily be unstable. This is again quite powerful; no further information is needed about the function $f$ or the pulse $\varphi$ that it supports. The topology is, in a sense, forcing the existence of a positive eigenvalue. 

Some remarks may be helpful here. First, the proof of the ``square" relies on the Maslov index and its topological properties, although the definition of $\theta$ given in \eqref{E:def-theta} is not directly used. Instead, the result is developed using the associated crossing form presented in \cite{RobbinSalamon93}. Second, a key step in the proof is proving a so-called monotonicity result. The path $\ell(x; \lambda) = \mathbb{E}^u_-(x; \lambda)$ is a path around the entire boundary of the square, if one considers either $x$ or $\lambda$ to be the path parameter on the appropriate sides, and hence a loop. After compactifying the domain, so that $x \in \mathbb{R}$ becomes $\tilde x \in [-1, 1]$, since the boundary of $[-1, 1] \times [0, \lambda_\infty] \subset \mathbb{R}^2$ is contractible, {\color{black} its image in the Lagrangian-Grassmannian is also contractible, and hence} the Maslov index of the loop $\ell(x; \lambda)$ must be zero. Showing there can be no intersections on the right side, where $\lambda = \lambda_\infty$ sufficiently large, or on the bottom, where $\tilde x = -1$, is not too difficult. One can then show that all crossings on the top (eigenvalues) must contribute in a negative way to the index, while on the left (conjugate points) they must contribute in a positive way; this is the monotonicity. {\color{black} Another way to view this monotonicity is in terms of the matrix $W$, defined in \eqref{E:def-theta}. In this setting, $W = W(x, \lambda)$, and this monotonicty result means that eigenvalues of $W$ must always pass through $-1$ in the same direction as $\lambda$ is varied, and always in the opposite direction as $x$ is varied.} Hence, the number of eigenvalues must equal the number of conjugate points. The fact that there must be at least one conjugate point when linearizing about a pulse comes from a symmetry argument that uses the reversibility of \eqref{E:main-rd} (the fact that it is invariant under the transformation $x \to -x$).

Not only does this result allow for the extension to the system case of the ``pulses must be unstable" result from Sturm-Liouville theory, but it also provides a more efficient way, in general, for detecting instabilities, provided one has the required symplectic structure, for example if $f = \nabla G$. To explain this, note that the Evans function, mentioned above, can be defined by $\mathcal{E}(\lambda) := \mathbb{E}^u_-(0; \lambda) \wedge  \mathbb{E}^s_+(0; \lambda): \mathbb{C} \to \mathbb{C}$. (The choice to look for intersections of the subspaces of decaying solutions at $x = 0$ is arbitrary; any point $x_0 \in \mathbb{R}$ could be chosen here.) Zeros of the Evans function correspond, with multiplicity, to eigenvalues. In general, to detect instabilities using the Evans function, one must prove that any unstable eigenvalues must lie in some compact ball and then compute the winding number of $\mathcal{E}$ around the boundary of this ball. On the other hand, to count conjugate points, one must only do analysis for a single value $\lambda = 0$. Thus, if one were to use validated numerics to produce a proof of (in)stability via such a detection procedure, the computation would be much faster using conjugate points than using the Evans function. This is the subject of current work.  

It is interesting to note that this connection between the Maslov index and stability, including the above demonstration of pulse instability, is not the only connection between topology and dynamic stability. It is also know that in some systems that support traveling waves, the wave can be constructed as the intersection of appropriate stable and unstable manifolds. This intersection typically occurs for a unique wavespeed, and the direction in which those manifolds cross as the wavespeed parameter is varied can be connected with $\mathcal{E}'(0)$ and hence the parity of the number of unstable eigenvalues; if this number is odd, there must be at least one, and the wave is unstable \cite{Jones95}. 

% Maybe email Chris or JB to ask for other examples of this?

So far, everything that has been discussed for \eqref{E:main-rd} has been restricted to the case of one spatial domain, $d = 1$. It turns out, however, that these ideas can also be expanded to cover the multidimensional case \cite{DengJones11, CoxJonesMarzuola15} %\cite{CoxJonesLatushkin16}
In this case, the eigenvalue problem takes the form
\begin{gather*}
\lambda v = \Delta v + \nabla^2G(\varphi(x)) v, \qquad x \in \Omega \subset \mathbb{R}^d \\
v|_{\partial \Omega} = 0.
\end{gather*}
To create the above theory in this setting, we need a notion of a conjugate point. This can be defined {\color{black} using a one-parameter family} of domains, $\{\Omega_s: s \in [0, 1], \Omega_1 = \Omega, \Omega_0 = \{x_0\}\}$, that shrinks the original domain down to a point \cite{Smale65}. One can then construct the path of subspaces
\begin{eqnarray*}
\ell(s; \lambda) &=& \left\{ \left( u, \frac{\partial u}{\partial n} \right)\Big|_{\partial \Omega_s}: u \in H^1(\Omega_s), \right. \\  
&& \qquad \left. \Delta u + V(x) u = \lambda u, \quad x \in \Omega_s\right\}
\end{eqnarray*}
determined by weak solutions on $\Omega_s$, but with no reference yet to the boundary data. {\color{black} By considering the Hilbert space
\[
\mathcal{H} = H^{1/2}(\partial \Omega) \times H^{-1/2}(\partial \Omega)
\]
and the symplectic form  $\omega((f_1, g_1), (f_2, g_2)) = \langle g_2, f_1 \rangle -  \langle g_1, f_2 \rangle$, where $\langle \cdot, \cdot \rangle$ denotes the dual pairing, once can show that both the path $\ell$ and the Dirichlet subpace
\[
\mathcal{D} = \left\{ \left( u, \frac{\partial u}{\partial n} \right) = \left( 0, \frac{\partial u}{\partial n} \right) \right\} \subset \mathcal{H},
\]
lie in the associated Fredholm-Lagrangian-Grassmanian, a generalization of the Lagrangian Grassmannian $\Lambda(n)$ to the infinite-dimensional setting. This Dirichlet subspace is now the fixed reference space, and a conjugate point is a value of $s$ such that $\ell(s; \lambda) \cap \mathcal{D} \neq \{0\}$. Note that the term ``Dirichlet subspace" in this context references the fact that this subspace corresponds to the zero Dirichlet boundary condition in the above eigenvalue problem. This perspective was pioneered in \cite{DengJones11} and allows for much of the above theory to work for the multi-dimensional eigenvalue problem, including the system case $v \in \mathbb{R}^n$ and a variety of boundary conditions other than Dirichlet.}

%The fixed reference space is again the Dirichlet space, 
%\[
%\mathcal{D} = \left\{ \left( u, \frac{\partial u}{\partial n} \right) = \left( 0, \frac{\partial u}{\partial n} \right) \right\},
%\]
%and so a conjugate point is a value of $s$ such that $\ell(s; \lambda) \cap \mathcal{D} \neq \{0\}$. Moreover, by considering the Hilbert space $\mathcal{H} = H^{1/2}(\partial \Omega) \times H^{-1/2}(\partial \Omega)$ and the symplectic form  $\omega((f_1, g_1), (f_2, g_2)) = \langle g_2, f_1 \rangle -  \langle g_1, f_2 \rangle$, where $\langle \cdot, \cdot, \rangle$ denotes the dual pairing, it can be shown that both the Dirichlet subspace and the path $\ell$ lie in the associated Fredholm-Lagrangian-Grassmanian, a generalization of the Lagrangian Grassmannian $\Lambda(n)$ to the infinite-dimensional setting. This perspective was pioneered in \cite{DengJones11} and allows for much of the above theory to work for the multi-dimensional eigenvalue problem, including the system case $v \in \mathbb{R}^n$ and a variety of boundary conditions other than Dirichlet. 

These multidimensional results are particularly exciting because most of the results related to nonlinear waves and coherent structures, not just their stability, apply only in one dimension. This is largely because many of the techniques rely on the perspective of spatial dynamics, which, for the most part, only applies to systems in one space dimension, or on cylindrical domains with a single distinguished spatial variable. Interestingly, the above procedure of using a shrinking family of domains, $\{\Omega_s\}$, suggests a way to develop spatial dynamics in higher dimensions. 

\section*{Spatial Dynamics}

{\color{black} In order to more precisely characterize what is meant by the term ``spatial dynamics," let's} recall the most basic setting in which spatial dynamics has been used, second order ordinary differential equations (ODEs) of the form $u_{xx} + F(u) = 0$. By writing this as the first order system
\[
u_x = v, \qquad v_x = -F(u),
\] 
one can study the behavior of solutions using techniques from dynamical systems, such as phase plane analysis and exponential dichotomies. Here the spatial domain is $\Omega = \mathbb{R}$, and the phase space of the spatial dynamical system is $\mathbb{R}^2$ (or $\mathbb{R}^{2n}$ if $u \in \mathbb{R}^n$). {\color{black} The above system is a spatial dynamical system, or equivalently it is the second order ODE viewed from the perspective of spatial dynamics, because in it the spatial variable $x$ is viewed as a time-like evolution variable, and techniques from the theory of dynamical systems can be used to study an equation that was not originally formulated as an evolutionary equation.}

On a cylindrical domain, $\Omega = \mathbb{R} \times \Omega'$ with $\Omega' \subset \mathbb{R}^{d-1}$ compact, the PDE $\Delta u + F(u) = 0$ can be written
\begin{equation}\label{E:infinite-sds}
u_x = v, \qquad v_x = -\Delta_{\Omega'}v - F(u),
\end{equation}
where $\Delta_{\Omega'}$ is the Laplacian on the cross section $\Omega'$. The phase space is now infinite-dimensional, for example $(u,v)(x) \in H^1(\Omega') \times L^2(\Omega')$ for each $x \in \mathbb{R}$, and so one must be more careful in analyzing the dynamics. This can be seen explicitly if $\Omega' = [0, 2\pi]$ with periodic boundary conditions, in which case the linear part of  \eqref{E:infinite-sds} coming from the Laplacian,
\[
\begin{pmatrix} 0 & 1 \\ -\partial_y^2 & 0 \end{pmatrix},
\]
has spectrum equal to the integers. This can be seen by using the Fourier expansion $u(x, y) = \sum_{k} \hat{u}_k(x)e^{\rmi k y}$, $v(x, y) = \sum_{k} \hat{v}_k(x)e^{\rmi k y}$, in which case $-\partial_y^2 \to k^2$ and the eigenvalues can be explicitly computed.
The fact that there are arbitrary large positive and negative eigenvalues means that, in general, solutions to \eqref{E:infinite-sds} will grow arbitrarily fast both forwards and backwards in time. In other words, the system \eqref{E:infinite-sds} is ill-posed. Nevertheless, 
applying techniques from dynamical systems to analyze the behavior of solutions is extremely useful. 

For example, in many cases one can construct an exponential dichotomy associated with the linear part of \eqref{E:infinite-sds}, and also construct stable and unstable (or possibly center-stable and center-unstable) manifolds associated with the nonlinear system. This allows for the analysis of subspaces, in the case of the dichotomy, or more generally manifolds of solutions that exist in forwards or backwards time, respectively. As a result, one can study bifurcations by looking at intersections of the relevant manifolds as system parameters are varied. One can also study stability, both at the spectral level using a generalization of the Evans function, at the linear level using pointwise Green's function estimates, and at the nonlinear level by combining these estimates with a representation of solutions to the full nonlinear equation, for example via Duhamel's formula. This infinite-dimensional spatial dynamics perspective began with the work of Kirchgassner \cite{Kirchgassner82}, and subsequent contributions include \cite{Mielke86, PeterhofSandstedeScheel97}. %Gallay93, %BeckSandstedeZumbrun10

The perspective of spatial dynamics has proven to be quite useful, and it has allowed for {\color{black} an extensive variety} of interesting and beautiful results to be obtained for PDEs on either one-dimensional or cylindrical domains. It has not, however, been utilized in multidimensional domains that do not have this cylindrical structure, and this is arguably the main reason why there are many fewer results available in higher space dimensions. The hope is that recent results, motivated by the above stability theory and which I will now describe, will change this.

Consider the PDE 
\begin{equation}\label{E:elliptic-pde}
\Delta u + F(x, u) = 0
\end{equation} 
with $x \in \Omega \subset \mathbb{R}^d$, and recall Smale's idea of shrinking the domain $\Omega$ via a one-parameter family $\{\Omega_s\}_{s\in[0,1]}$. Suppose that this family is parameterized by a family of diffeomorphisms $\psi_s: \Omega \to \Omega_s$. This allows for a nice definition of the boundary data on $\partial \Omega_s$:
\[
f(s; y) = u(\psi_s(y)), \qquad g(s; y) = \frac{\partial u}{\partial n}(\psi_s(y)), 
\]
for $s \in [0, 1]$ and $y \in \partial \Omega$. This is convenient because, even though $(f,g)$ can be interpreted as the boundary data on $\Omega_s$, the independent variable $y$ lives in the $s$-independent domain $\partial\Omega$. One can then, at least formally, compute an evolution equation of the form
\begin{equation}\label{E:ses}
\frac{d}{ds}(f, g) = \mathcal{F}(f, g),
\end{equation}
where the possibly nonlinear function $\mathcal{F}$ is defined in terms of the function $F$ appearing in \eqref{E:elliptic-pde} and the tangential parts of the gradient and divergence operators on $\partial \Omega_s$. One can also, again at least formally, relate a solution $(f,g)$ of \eqref{E:ses} to the solution $u$ of \eqref{E:elliptic-pde} {\color{black} by noting that $(f,g)$ is just the function $u$ and its normal derivative evaluated on the boundary on the domain $\Omega_s$; in other words, $(f,g)$ is just the trace of $(u, \partial u/\partial n)$ evaluated on $\partial \Omega_s$.}
%via the trace operator, 
%\[
%\mathrm{Tr}_s u = (f(s), g(s)).
%\]
This has been made rigorous in \cite{BeckCoxJones19a}, where it was shown that, in an appropriate sense, a weak solution $u$ of the elliptic PDE \eqref{E:elliptic-pde} leads to a solution $(f, g)$ of the spatial dynamical system \eqref{E:ses}, and vice versa. 

The function $\mathcal{F}$ is indeed quite complicated, and the relation between $u$ and $(f, g)$ is rather technical. However, for at least some domains $\Omega$, the result seems to be sufficiently concrete so as to be {\color{black} readily applicable}. For example, if the domain is radial or all of $\mathbb{R}^d$, one can choose to shrink the domain using spheres: $\Omega_s = \{x \in \mathbb{R}^d: |x| < s\}$. This greatly simplifies the function $\mathcal{F}$ and, using the fact that in terms of generalized polar coordinates $\Delta = \partial_r^2 + (n-1)r^{-1} \partial_r+ r^{-2}\Delta_{\mathbb{S}^{d-1}}$, one ends up with the spatial dynamical system 
\begin{eqnarray*}
\frac{d}{ds} \begin{pmatrix} f \\ g \end{pmatrix} &=& \begin{pmatrix} 0 & 1 \\ - s^{-2} \Delta_{\mathbb{S}^{d-1}} & -(d-1)s^{-1}\end{pmatrix} \begin{pmatrix} f \\ g \end{pmatrix} \\
&& \qquad \qquad + \begin{pmatrix} 0 \\ - F(\theta, s, f) \end{pmatrix}.
\end{eqnarray*}
It has been shown that the linear part of this system, after a suitable rescaling of time $s = e^\tau$ and for $d \geq 3$, admits an exponential dichotomy \cite{BeckCoxJones19b}. (The case $d = 2$ is slightly more complicated, due to the existence of the harmonic function $\log r$, but it could be similarly interpreted by allowing the dichotomy to contain center directions.) Moreover, when $d=3$ the dichotomy can be written down explicitly in terms of the {\color{black}spherical harmonics}. This allows one to potentially study solutions to the original elliptic PDE that are not necessarily radially symmetric, thus providing the removal of a restriction that has been imposed on most results (at least in the spatial dynamics context) to date. Thus, the perspective of spatial dynamics seems quite promising as a method for studying multidimensional nonlinear waves and coherent structures. 

%{\color{black} I could add in more details here about the rescaling or the exponential dichotomy. If I have space I could also add more figures or more not-in-line equations.}

\section*{Future Directions}

The theory discussed above has the potential to have a great impact, particularly for problems in multiple spatial dimensions. Many of the existing results are valid only for one-dimensional domains, or for cylindrical domains. The above results represent new techniques that are not bound by this restriction, and thus allow for the analysis not only of stability but also of a variety of aspects of the behavior of solutions to PDEs in multi-dimensional spatial domains, {\color{black} such as their existence and bifurcation}. 

In the last ten years or so there have been many results regarding the theory discussed above. Arguably the only downside so far is the relative lack of applications: examples of solutions, in any space dimension, whose stability is determined using the conjugate point method described above and instances of using the spatial dynamical system \eqref{E:ses} to analyze multidimensional nonlinear waves.

Regarding the former, there are three existing examples, at least where the Evans function cannot also be used to determine stability. The most broadly applicable is the pulse instability result in reaction-diffusion systems with gradient nonlinearity, described above. The other two examples pertain to specific PDEs, with the first being the instability result of \cite{Jones88} for a standing wave in a nonlinear Schr\"odinger-type equation, which really began this whole program, and the second being the instability result of \cite{ChenHu14} for a standing pulse in the FitzHughÐNagumo equation, with diffusion in both variables. The development of the spatial dynamical system \eqref{E:ses} and its relation to the elliptic PDE \eqref{E:elliptic-pde} is extremely new, and so some time is needed for its utility to be fully explored. Now that a solid foundational theory is in place, the hope is that many more applications will emerge. This is an area of active, ongoing work.   

%in There is a related result where the Maslov index is used to prove that the standing pulse solution of the FitzHughÐNagumo equation, with diffusion in both variables, is unstable, provided certain assumptions on the model parameters are satisfied \cite{ChenHu14}. Unlike the rather general pulse instability result, however, the analysis of \cite{ChenHu14} is particular to the FitzHugh-Nagumo model and relies heavily on the activator-inhibitor structure in evaluating the Maslov index. Similarly, the instability result of \cite{Jones88}, mentioned above, is specific to the nonlinear Schr\"odinger equation analyzed there. 

%{\color{black} Read Gene's notices articles to see the level they're written for - 2nd year grad students? Take a look at all the paper's I'm citing, particularly the classical ones I don't really know well. Cut citations down to 20!!}

%%%%%%%%%%%%%%%%%%%%%%%%%%%%%%%%%%%%%%%%%%%%%%%%%%%%%%

\bibliography{notices.bib} % Limit of 20 references.

%The AMS uses MR style for references. If you are using a BibTeX style
%bibliography file (.bib), these will be formatted automatically by this
%document. Otherwise, follow the guides below. \textbf{Please also be sure
%to include MR numbers, if they exist, for all references listed.} 
%If you are cutting and pasting your references from MathSciNet, these are
%included automatically. Otherwise, please add the MR number as in this one
%(from a nice article on Galois theory for quasi-fields \cite{MR0001219}):

%If you don't have an MR number, the DOI and URL fields are also very useful for
%the Notices staff. The URL field can contain any link which will take the
%user to the paper. 

\end{document}